\documentclass[a4paper,12pt]{article}
\usepackage{amsmath}
\usepackage{amssymb}
\usepackage[dvips]{graphicx}

\newtheorem{theorem}{Theorem}[section]

\newtheorem{proposition}[theorem]{Proposition}

\newtheorem{definitionhead}[theorem]{Definition}
\newenvironment{definition}{\begin{definitionhead}%
\sl}{\end{definitionhead}}

\newcommand\Proof{\noindent{\bf Proof. }}
\newcommand\Endproof{\nopagebreak\strut%
\nopagebreak\hfill\nopagebreak$\Box$\medbreak}
\newcommand\Noproof{\nopagebreak\strut%
\nopagebreak\hfill\nopagebreak$\Box$\par}

\def\Fol{\operatorname{Fol}}
\def\Card{\operatorname{Card}}

\begin{document}

\title{Describing the set of words generated by interval exchange transformation}
\author{A.Ya.~Belov \and A.L.~Chernyat'ev}
\date{}

\maketitle

\newpage
\pagestyle{plain} \tableofcontents
\newpage

\section{Introduction}
Methods of symbolic dynamics are rather useful in the study of combinatorial
properties of words, investigation of problems of number theory and theory of
dynamical systems. Let $M$ be a compact metric space, $U\subset M$ be its open
subspace, $f:M\to M$ be a homeomorphism of the compact into itself, and $x\in
M$ be an initial point. With the sequence of iterations, one can associate an
infinite binary word

$$
w_n=\left\{
\begin{array}{rcl}
   a,\ f^{(n)}(x_0)\in U\\
   b,\ f^{(n)}(x_0)\not\in U\\
\end{array}\right.
$$
which is called the {\it evolution} of point $x_0$. Symbolic dynamics
investigates the interrelation between the properties of the dynamical system
$(M,f)$ and the combinatorial properties of the word $W_n$. For words over
alphabets which comprise more symbols, several characteristic sets should be
considered: $U_1,\ldots ,U_n$.

By the {\it direct problem} of symbolic dynamics we mean the study of
combinatorial properties of the words generated by a given dynamical system;
the {\it inverse problem} of symbolic dynamics refers to the investigation of
the properties of the dynamical system, i.e., the properties of the compact set
$M$ and transform $f$, by the combinatorial properties of the word $W$.

Inverse problems of symbolic dynamics related to the unipotent transformation
of a torus were studied in paper~\cite{BK}.

Problems (both direct and inverse) related to the rotation of a circle bring
about a class of words which are called {\it Sturmian words}. Sturmian words
are infinite words over a binary alphabet which contain exactly $n+1$ different
subwords (factors) of length $n$ for any $n\geq 1$. The following classical
result is widely known.

\begin{theorem}[Equivalence theorem (\cite{MorseHedlund},\cite{L1}).]\label{ThEquivClassik}
Let $W$ be an infinite recurrent word over the binary alphabet $A=\{a,b\}$. The
following conditions are equivalent:

\begin{enumerate}
\item The word $W$ is a Sturmian word, i.e., for any $n\geq 1$,
the number of different subwords of length $n$ that occur in $W$ is equal to $T_n(W)=n+1$.

\item The word is not periodic and is {\em balanced}, i.e., any two subwords $u,v\subset W$
of the same length satisfy the inequality $||v|_a-|u|_a|\leq 1$, where $|w|_a$
denotes the number of occurrences of symbol $a$ in the word $w$.

\item The word $W=(w_n)$ is a mechanical word with irrational $\alpha$,
which means that there exist an irrational $\alpha$, $x_0 \in [0,1]$,
and interval $U\subset \mathbb{S}^1$, $|U|=\alpha$,
such that the following condition holds:
$$
w_n=\left\{
\begin{array}{rcl}
   a,&{T_{\alpha}}^n(x_0)\in U\\
   b,&{T_{\alpha}}^n(x_0)\not\in U\\
\end{array}\right.
$$
\end{enumerate}
\end{theorem}

There are several different ways of generalizing Sturmian words.

First, one can consider {\it balanced} words over an arbitrary alphabet.
Balanced nonperiodic words over an $n$-letter alphabet were studied in
paper~\cite{GR} and later in~\cite{H}. In papers~\cite{BC} and \cite{C1}, a
dynamical system that generates an arbitrary nonperiodic balanced word was
constructed.

Second, generalization may be formulated in terms of the {\it complexity
function}. Complexity function $T_W(n)$ presents the number of different
subwords of length $n$ in the word $W$. Sturmian words satisfy the relation
$T_W(n+1)-T_W(n)=1$ for any $n\geq 1$. Natural generalizations of Sturmian
words are words with minimal growth, i.e., words over a finite alphabet that
satisfy the relation $T_W(n+1)-T_W(n)=1$ for any $n\geq k$, where $k$ is a
positive integer. Such words were described in terms of rotation of a circle in
paper~\cite{C2}. Note also that words whose growth function satisfies the
relation $\lim_{n \to \infty} T(n)/n = 1$ were studied in paper~\cite{Ab}.

Words with complexity function $T_W(n)=2n+1$ were studied by P.~Arnoux and
G.~Rauzy (\cite{AR,Ra,Ra1}), words with growth function $T_W(n)=2n+1$ were
investigated by G.~Rote~\cite{R}. Consideration of the general case of words
with linear complexity function involves the study of words generated by
interval exchange transformations. The problem of describtion of such words was
posed by Rauzy \cite{Ra1}.
Words with linear growth of the number of subwords were also studied
by the school of V.~Berth\'e, S.~Ferenczi, and Luca Q.~Zamboni
(\cite{FZ}, \cite{BFZ}). They also investigate combinatorial sequences reated
with interval exchange transformations. Paper~\cite{FHZ} contains description of words generated by
three-interval exchange transformations; paper~\cite{FZ} contains
description of words generated by symmetric interval exchange
transformations (such transformations are closely related to
multi-dimensional continued fractions, and this relation looks extremely intersting).
More precisely,   they
describe a combinatorial algorithm for
generating the symbolic sequences which code the orbits of points
under an interval exchange transformation on $k$ intervals, using
the “symmetric” permutation $i\to k-i+1$ (\cite{FZ}).


More general result was obtained in the work \cite{FZ3}. In this paper
give a complete characterization of those sequences of subword
complexity $(k - 1)n + 1$ which are natural codings of orbits of
$k$-interval exchange transformations, (or, equivaently, interval exchange transformation, satisfying i.d.o.c. condition).
\footnote{This result wich is close to ours, was noticed to us by L.Zamboni, in order to mention
that and make some other correstions we replaced our paper to new version}.

Interval exchange transformations $T$ satisfies {\it the infinite
distinct orbit condition} (or {\it i.d.o.c.} for short) if the $k-1$
negative trajectories $\{T^{-n}(x_i)\}_{n\geq 0}$ ,$(1\leq i \leq
k)$, of the discontinuities of $T$ are infinite disjoint sets. The
main result of paper \cite{FZ3} (wich is independent of our result) is following:

\begin{theorem}
A minimal sequence $W$ is the natural coding of a
$k$-in\-ter\-val exchange transformation, defined by permutations
$(\pi_0,\pi_1)$ such that ${\pi_0}^{-1}(\{1,\dots,j\})\neq
{\pi_1}^{-1}(\{1,\dots,j\})$ for every $1\leq  j \leq k-1$, and
satisfying the i.d.o.c. condition, if and only if the words of
length one occurring in $W$ are $F_1 = \{1,\dots,k\}$ and it
satisfies the following conditions:
\begin{enumerate}
\item If $w$ is any word occurring in $W$, $A(w)$, (resp. $D(w)$), the set of
all letters $x$ such that $xw$, (resp. $wx$), occurs in $W$, is an
interval for the order of $\pi_1$, resp. $\pi_0$,

\item If $x \in A(w)$, $y \in A(w)$, $x \leq y$ for the order of $\pi_1$, $z \in D(xw)$,
$t\in D(yw)$, then $z \leq t$ for the order of $\pi_0$,

\item If $x \in A(w)$ and $y \in A(w)$ are consecutive in the order of $\pi_1$,
$D(xw) \cap D(yw)$ is a singleton.
\end{enumerate}
\end{theorem}

In this paper we study words generated by general piecewise-continuous transformation of the interval.
Further we prove equivalence set words generated by piecewise-continuous transformation
and words generated by interval exchange transformation.
This method get capability of descriptions of the words generated by
arbitrary interval exchange transformation.

This work is targeted to the following

\medskip
{\bf Inverse problem}: {\it Which
conditions should be imposed on a uniformly recurrent word $W$ in order that it
be generated by a dynamical system of the form $(I,T,U_1,\ldots,U_k)$, where
$I$ is the unit interval and $T$ is the interval exchange transformation?
}
\medskip

The answer to this question is given in terms of the evolution of
the {\it labeled Rauzy graphs} of the word $W$. The {\it Rauzy
graph} of order $k$ (the $k$-graph) of the word $W$ is the directed
graph whose vertices biuniquely correspond to the factors of length
$k$ of the word $W$ and there exists an arc from vertex $A$ to
vertex $B$ if and only if $W$ has a factor of length $k+1$ such that
its first $k$ letters make the subword that corresponds to $A$ and
the last $k$ symbols make the subword that corresponds to $B$. By
the {\it follower} of the directed $k$-graph $G$ we call the
directed graph $\Fol(G)$ constructed as follows: the vertices of
graph $\Fol(G)$ are in one-to-one correspondence with the arcs of
graph $G$ and there exists an arc from vertex $A$ to vertex $B$ if
and only if the head of the arc $A$ in the graph $G$ is at the notch
end of $B$. The $(k+1)$-graph is a subgraph of the follower of the
$k$-graph; it results from the latter by removing some arcs.
Vertices which are tails of (or heads of) at least two arcs
correspond to {\it special factors} (see Section~2); vertices which
are heads and tails of more than one arc correspond to {\it
bispecial factors}. The sequence of the Rauzy $k$-graphs constitutes
the {\it evolution} of the Rauzy graphs of the word $W$. The Rauzy
graph is said to be {\it labeled} if its arcs are assigned letters
$l$ and $r$ and some of its vertices (perhaps, none of them) are
assigned symbol ``-''.

The {\it follower} of the labeled Rauzy graph is the directed graph
which is the follower of the latter (considered a Rauzy graph with
the labeling neglected) and whose arcs are labeled according to the
following rule:

\begin{enumerate}

\item Arcs that enter a branching vertex should be labeled by the
same symbols as the arcs that enter any left successor of this
vertex;

\item Arcs that go out of a branching vertex should be labeled by
the same symbols as the arcs that go out of any right successor of
this vertex;

\item If a vertex is labeled by symbol ``--'', then all its right
successors should also be labeled by symbol ``--''.

\end{enumerate}

In terms of Rauzy labeled graphs we define the {\it
asymptotically correct} evolution of Rauzy graphs, i.e., we
introduce rules of passing from $k$-graphs to $(k+1)$-graphs.
Namely, the evolution is said to be {\it correct} if, for all
$k\geq 1$, the following conditions hold when passing from the
$k$-graph $G_k$ to the $(k+1)$-graph $G_{k+1}$ :

\begin{enumerate}
\item The degree of any vertex is at most $2$, i.e., it is incident to at most two
incoming and outgoing arcs;

\item If the graph contains no vertices corresponding to bispecial factors,
then $G_{n+1}$ coincides with the follower $D(G_n)$;

\item If the vertex that corresponds to a bispecial factor is not
labeled by symbol ``--'', then the arcs that correspond to
forbidden words are chosen among the pairs $lr$ and $rl$;

\item If the vertex is labeled by symbol ``--'', then the arcs to
be deleted should be chosen among the pairs $ll$ or $rr$.
\end{enumerate}

The evolution is said to be {\it asymptotically correct} if this
condition is valid for all $k$ beginning with a certain $k=K$. The
{\it oriented} evolution of the graphs means that there are no
vertices labeled by symbol ``--''. The main result of this work
consists in the description of infinite words generated by
interval exchange transformations (and answers a Rouzy question \cite{Ra1}):

\medskip
{\bf Main theorem.}\ {\it A uniformly recurrent word $W$
 \begin{enumerate}
\item is generated by an interval exchange transformation if and
only if the word is provided with the asymptotically correct
evolution of the labeled Rauzy graphs;

\item is generated by an or\-i\-en\-ta\-ti\-on--pres\-er\-ving interval exchange
transformation if and only if the word is provided with the
asymptotically correct oriented evolution of the labeled Rauzy
graphs.
\end{enumerate}
}
\medskip

We have no restriction on the endpoint orbit. In special case of
asymptotical subword growth of $T_W(n)=n+\mbox{const}$ for all
$n>n_0$ we get an generalization of theorem \ref{ThEquivClassik},
i.e. description of all u.r. words with such growth property. The
description of all (not nesesary u.r.) superwords such that $T_W(n)=n+\mbox{const}$
for all
$n>n_0$ see in \cite{C2}.
 Note also that in all previous studies which is known for us interval
 exchange transformations are defined to be orientation preserving.

The paper is organized as follows: in Section 2 we formulate the main definitions and facts
about uniformly recurrent words, Rauzy graphs, and words generated by dynamical systems.
In Section 3 we prove a theorem about necessary conditions for a word to be generated by
interval exchange transformation. The next two sections are devoted
to the proof of the sufficiency of these conditions. In Section 4
we prove that these conditions are sufficient for the word
to be generated by a piecewise-continuous interval transformation. Finally,
in Section 5 we prove that the sets of uniformly recurrent words generated by
piecewise-continuous interval transformations and by the interval exchange transformation
are equivalent.

\section{Main constructions and definitions}

\subsection{Complexity function, special factors, and uniformly recurrent words}
In this section we define the basic notions of combinatorics of words.
By $L$ we denote a finite alphabet, i.e., a nonempty set of elements (symbols).
We use the notation $A^+$ for the set of all finite sequences of
symbols or {\it words}.

A finite word can always be uniquely represented in the form
$w = w_1 \cdots w_n$, where $w_i \in A$, $1\leq i \leq n$. The number $n$
is called the {\it length} of word $w$; it is denoted by $|w|$.

The set $A^+$ of all finite words over $A$ is a simple semigroup
with concatenation as semigroup operation.

If element $\Lambda$ (the empty word) is included in the set of
words, then this is actually the free monoid $A^*$ over $A$. By
definition the length of the empty word is $|\Lambda|=0$.

A word $u$ is a {\it subword} (or {\it factor}) of a word $w$ if there exist
words $p,q\in A^+$ such that $w = puq$.

Denote the set of all factors (both finite and infinite) of a word $W$ by $F(W)$.
Two infinite words $W$ and $V$ over alphabet $A$ are said to be {\it
equivalent} if $F(W)=F(V)$.

We say that symbol $a\in A$ is a {\it left} (accordingly, {\it right}) {\it extension} of factor $v$ if $av$
(accordingly, $va$) belongs to $F(W)$. A subword $v$ is called a {\it left} (accordingly, {\it right}) {\it special
factor} if it possesses at least two left (right) extensions. A subword $v$ is said to be {\it bispecial}
if it is both a left and right special factor at the same time.
The number of different left (right) extensions of a subword
is called the {\it left (right) valence} of this subword.

A word $W$ is said to be {\it recurrent} if each its factor
occurs in it infinitely many times (in the case of a doubly-infinite
word, each factor occurs infinitely many times in both directions).
A word $W$ is said to be {\it uniformly recurrent} or ({\it u.r word})
if it is recurrent and, for each its factor $v$, there exists a positive integer
$N(v)$ such that, for any subword $u$ of length at least $N(v)$ of the word $W$,
factor $v$ occurs in $u$ as a subword.

Below we formulate several theorems about u.r words, which will be needed later.
The proof of these theorems can be found in monograph~\cite{BBL}.

\begin{theorem}
The following two properties of an infinite word $W$
are equivalent:

a) For any $k$ there exists $N(k)$ such that any segment of length
$k$ of the word $W$ occurs in any segment of length $N(k)$ of the word $W$;

b) If all finite factors of a word $V$ are at the same time finite factors of a word $W$,
then all finite factors of the word $W$ are also finite factors of the word $V$.
\end{theorem}

\begin{theorem}\label{ThUnifRec}
Let $W$ be an infinite word. Then there exists a uniformly recurrent word $\widehat{W}$
all of whose factors are factors of $W$. \Noproof
\end{theorem}

One can consider the action of the shift operator $\tau$ on the set
of infinite words. {\it The Hamming distance} between words $W_{1}$
and $W_{2}$ is the quantity $d(W_{1},W_{2}) =\sum_{n\in {\mathbb
Z}}\lambda _{n}2^{-|n|}$, where $\lambda _{n}= 0$ if symbols at the
$n$-th positions of the words are the same and $\lambda _{n}= 1$,
otherwise.

An {\it invariant subset} is a subset of the set of all infinite
words which is invariant under the action of $\tau$. A {\it
minimal closed invariant set}, or briefly, m.c.i.s, is a closed
(with respect to the Hamming metric introduced above) invariant
subset which is nonempty and contains no closed invariant subsets
except for itself and the empty subset.

\begin{theorem}[Properties of closed invariant sets]\label{RecEq}
The following properties of a word $W$ are equivalent:
\begin{enumerate}
\item $W$ is a uniformly recurrent word;

\item The closed orbit of $W$ is minimal and is a m.c.i.s.
\end{enumerate}
\end{theorem}

\begin{theorem}\label{RecCont}
Let $W$ be a uniformly recurrent nonperiodic infinite word.
Then
\begin{enumerate}
\item All the words that are equivalent to $W$ are u.r. words;
the set of such words in uncountable;

\item There exist distinct u.r. words $W_1 \neq W_2$ which are equivalent to the given word
and can be written as $W_1=U V_1$, $W_2=U V_2$, where $U$ is a left-infinite word
and $V_1\neq V_2$ are right-infinite words.
\end{enumerate}
\end{theorem}

\subsection{Rauzy graphs}
It is convenient to describe a word $W$ using the {\it subword
graphs} or (Rauzy graphs). They were introduced by Rauzy~\cite{Ra}
in the following way: the {\it $k$-graph} of the word $W$
is the directed graph whose vertices biuniquely correspond
to subwords of length $k$ of the word $W$ and there is an arc from vertex
$A$ to vertex $B$ if $W$ contains a subword of length $k+1$
such that the first $k$ symbols of it make a subword that corresponds to $A$
and the last $k$ symbols make a subword that corresponds to $B$.
Thus, the arcs of the $k$-graph are in one-to-one correspondence with
the ($k+1$)-factors of the word $W$.

It is clear that, in the $k$-graph $G$ of the word $W$, vertices
which are tails (accordingly, heads) of more than one arc correspond to right special words.
Such vertices will be called crotches. Graph $G$ is said to be {\it strongly connected}
if it contains a directed path from any vertex to any other one.

The {\it follower} of the directed graph $G$ is the directed graph
$\Fol(G)$ constructed in the following way:
the vertices of the graph $G$ biuniquely correspond to the arcs of the $G$
and there is an arc from vertex $A$ to vertex $B$ if the head
of the arc $A$ in the graph $G$ is at the notch end of $B$.

The connectivity of the Rauzy graphs is naturally related to the
recurrence of the corresponding word. Namely, the following assertion is valid.

\begin{proposition}\label{reccur}
Let $W$ be a (semi)infinite word. The following conditions are equivalent:

\begin{enumerate}
\item The word $W$ is recurrent;

\item For any $k$ the corresponding $k$-graph of the word $W$ is strongly connected;

\item Any factor of $W$ occurs at least twice;

\item Any factor can be extended to the left.
\end{enumerate}
\end{proposition}

Let us introduce the notion of the {\it labeled Rauzy graph}. A
Rauzy graph is said to be {\it labeled} if

\begin{enumerate}
\item the arcs of any crotch are assigned symbols $l$ (``left'') and $r$ (``right'');

\item some vertices are assigned symbol ``--''.
\end{enumerate}

The {\it follower} of the labeled Rauzy graph is the directed graph
which is the follower of the latter (considered a Rauzy graph with
the labeling neglected) and whose arcs are labeled according the
following rule:
\begin{enumerate}

\item Arcs that enter a crotch should be labeled by the same
symbols as the arcs that enter any left successor of this vertex;

\item Arcs that go out of a crotch should be labeled by the same
symbols as the arcs that go out of any right successor of this
vertex;

\item If a vertex is labeled by symbol ``--'', then all its right
successors should also be labeled by symbol ``--''.

\end{enumerate}

\subsection{Words generated by dynamical systems}

Let $M$ be a compact metric space, $U\subset M$ be its open subset, $f:M\to M$
ba a homeomorphism of the compact space into itself, and $x\in M$ be an initial point.

With the sequence of iterations, one can associate an
infinite binary word

$$
w_n=\left\{
\begin{array}{rcl}
   a,\ f^{(n)}(x_0)\in U\\
   b,\ f^{(n)}(x_0)\not\in U\\
\end{array}\right.
$$
which is called the {\it evolution} of point $x_0$. Symbolic dynamics
investigates the interrelation between the properties of the dynamical system
$(M,f)$ and the combinatorial properties of the word $W_n$.

For words over alphabets which comprise more symbols,
several characteristic sets should be considered: $U_1,\ldots ,U_n$.

Note that the evolution of the point is correctly defined
only when the trajectory of the point does not pass through the boundary of
the characteristic sets $\partial U_1,\partial U_2,\ldots$.

In order to consider the trajectory of an arbitrary point,
let us introduce the notion of {\it essential evolution}.

\begin{definition}
A finite word $v^f$ is called {\it essential finite evolution} of
point $x^{\ast}$ if any neighborhood of point $x^{\ast}$ contains
an open set $V$ such that any point $x\in V$ possesses the
evolution $v^f$. An infinite word $W$ is called {\it essential
evolution} of point $x^{\ast}$ if any its initial subword is an
essential finite evolution of point $x^{\ast}$.
\end{definition}

When there is no risk of ambiguity, we say {\it evolution of the point}
meaning the essential evolution.
Note that a point can have several essential evolutions.

\begin{proposition}[\cite{BK}] Let $V$ be a finite word. Then the set of points with
the finite essential evolution $V$ is closed. A similar assertion is true
for an infinite word $W$.
\end{proposition}

\subsection{Correspondence between words and partitions of a set}\label{Sootv:2}
Now let us consider the correspondence between words and subsets
of $M$. It follows from the construction that, if the initial
point belongs to the set $U_i$, then its evolution begins with
symbol $a_i$. Consider the images of the sets $U_i$ under the
mappings $f^{(-1)}, f^{(-2)}, \ldots$, $n\in \mathbb{N}$.

It is clear that, if the point belongs to the set
$$
   f^{(-n)}(U_{i_n}) \cap f^{-(n-1)}(U_{i_{n-1}}) \cap \ldots \cap f^{(-1)}(U_{i_1}) \cap U_{i_0},
$$
then the evolution begins with the word $a_{i_0}a_{i_1}\cdots a_{i_n}$.

Accordingly, the number of different essential evolutions of length
$n+1$ is equal to the number of partitions of the set $M$ into nonempty subsets
by the boundaries of the sets $\partial U_i$ and their images under the mappings
$f^{(-1)}, f^{(-2)}, \ldots ,f^{(-n)}$.

{\bf Remark.} The number of finite essential evolutions is
directly related to the topological dimension of the set $M$. For
instance, it is clear that, if $і$ is homeomorphic to a segment or
a circle, then one point can belong to the boundary of at most two
open subsets of $M$ and, accordingly, can have only two essential
evolutions. If $M$ is homeomorphic to a part of the plane
$\mathbb{R}^2$ then there can be arbitrarily many essential
evolutions.


\subsection{Interval exchange transformations}\label{PerOtr}

Interval exchange transformation is a natural generalization of
the shift of a circle: in the case of the partition of a circle
into arcs of length $\alpha$ and $1-\alpha$ and a shift of
$\alpha$, this transformation coincides with the two-interval
exchange transformation.

In addition, interval exchange transformation is rather important in
ergodic theory, theory of dynamical systems, and number theory.

Let us consider the general case:

Suppose that the closed interval $[0,1]$ is partitioned into
half-open intervals of lengths $\lambda_1,\lambda_2,\ldots
,\lambda_n$ and $\sigma\in S_n$ is a permutation of the set
$\{1,2,\ldots,n\}$.

The intervals of the partition can be represented in terms of
the lengths given above:

$$
X_i=\left [ \sum_{i<j} \lambda_j,\sum_{i\leq j} \lambda_j\right).
$$

The interval exchange transformation rearranges the intervals $(X_1,X_2,\ldots,X_n)$
of the partition; as a result, we obtain a new partition

$$
(X_\sigma(1),X_\sigma(2), \ldots ,X_\sigma(n)).
$$

In the orientation-preserving case, transformation $T$ associates
each point $x\in X_i$ :
$$
T(x)=x+a_i,
$$
where
$$
a_i=\sum_{k<\sigma^{-1}(i)}\lambda_{\sigma(k)}-\sum_{k<i}\lambda_k.
$$

If the transformation inverts an interval, then, in addition, all points
are symmetrically reflected with respect to the midpoint of this segment.

\begin{definition}
Interval exchange transformation $T$ is said to be {\it regular}
if, for any point $a_i$, where $X_i=[a_i,a_{i+1})$, we have $T^n(a_i)\neq a_j$.
\end{definition}

The result formulated below is rather important (see~\cite{FHZ}):

\begin{theorem}
An interval exchange transformation is regular if and only if
the trajectory of any point is dense everywhere in $[0,1]$.
\end{theorem}

Properties of words generated by interval exchange transformations
are investigated using the same methods as in the case of the shift of the unit circle.
Here, the main approach consists in considering the negative orbits
of the ends of the exchanged intervals

$$
0 = a_1 < a_2 < \ldots < a_{k+1} = 1,
$$

where

$X_i = [a_i, a_{i+1})$, ($i \in \{1, \ldots , k\}$).

Denote the set of ends of the exchanged intervals $\{a_i| 1 \leq i \leq k
+ 1\}$ by $X^{1}$. A word $w$ of length $n$ is a subword
of the evolution of point $x$, i.e., the infinite word $U(x)$, if and only if
there exists an interval $I_w \subset[0, 1]$ and point $y \in I_w$
such that the word $w$ is the concatenation of the symbols:

$$
\mathcal{I}(x)\mathcal{I}(T(x))\cdots \mathcal{I}(T^{n-1}(x)) = w,
$$
where $\mathcal{I}(x)=a_i \in A $ if and only if $x\in X_i$.

 \begin{proposition}
 Let $T$ be a regular $k$-interval exchange transformation.
Then the evolution $U(x)$ of any point $x$ has the complexity function
$T_{U(x)}(n)=n(k-1)+1$ for any $n\in\mathbb{N}$.
\end{proposition}

\section{Necessary conditions for a word to be generated
by an interval exchange transformation}

At the first stage we formulate necessary conditions for a word to be
generated by an interval exchange transformation.
Let the word $W$ be the evolution of point $x\in [0,1]$
for a $k$-interval exchange transformation and
characteristic sets $U_1,U_2,\ldots,U_n$. Each characteristic set $U_i$
is a union of several disjoint open or half-open intervals.

As was already shown in Section~\ref{Sootv:2}, subwords of length $k$
are in one-to-one correspondence with the $k$-partitions of the characteristic sets.
Since a boundary point of a one-dimensional set can belong to
the boundary of only two sets, a $k$-subword of the word
$W$ can have at most two extensions. We obtain the first necessary condition
for a word to be generated by interval exchange transformation:

\begin{proposition}
Let the word $W$ be generated by interval exchange transformation.
Then, for a certain $N$, all special subwords of length at least $N$
should have valence $2$.
\end{proposition}

This condition is similar to the condition that, starting from a certain $N$,
all $k$-graphs of the word $W$ ($k\geq N$) should have all incoming and outgoing
crotches of degree $2$.

Now let us derive conditions in terms of Rauzy graphs. Suppose that the recurrent word
$W$ has the growth function $F_W(n)=Kn+L$ for $n>N$ and is generated
by interval exchange transformation. Consider the evolution of Rauzy $k$-graphs
of the word $W$ starting with $k\geq N+1$. As has already been demonstrated, all incoming
and outgoing crotches have degree $2$; therefore, all $k$-graphs have exactly
$K$ incoming and $K$ outgoing crotches.

In the minimal case described in the previous section, where
$F_W(n)=n+L$, Rauzy graphs have exactly one incoming and one outgoing
crotch. If the incoming crotch coincides with the outgoing one, then
the choice of the arc to be deleted from the follower $D(G)$
is uniquely determined by the condition of strong connectivity of the graph.

Let us thoroughly investigate a more interesting case, where graphs of words contain
more than one crotch. There are several possible situations that should be considered
when passing from graph $G_n$ to $G_{n+1}$:

\begin{enumerate}
\item Graph $G_n$ contains no linked cycles (i.e., there are no incoming crotches
that are at the same time outgoing crotches). In this case, graph
$G_{n+1}$ coincides with the follower $D(G_n)$.

\item Graph $G_n$ contains one crotch which is at the same time
an incoming and outgoing crotch. In this case, the follower graph $D(G_n)$
has three crotches, because one crotch has been cloned.
Therefore, the graph $G_{n+1}$ is obtained from the follower $D(G_n)$
by means of deleting one arc which corresponds to the minimal
non-occurring word.

\item The graph contains two or more crotches which are at the same time
incoming and outgoing crotches. Then the graph $G_{n+1}$ is obtained from $D(G_{n})$
by means of deleting two or more arcs which correspond to the minimal
non-occurring words.

\end{enumerate}

Since the word $W$ is recurrent, it follows from Proposition~\ref{reccur}
that, as the arcs are deleted, the graph should remain strongly connected,
i.e., it should contain a directed path from any vertex to any other one.

Consider the second case in more detail. Suppose that $G_k$
contains one double crotch. This means that $W$ contains exactly
one bispecial subword $w$ of length $k$. Hence, there exist
$a_i,a_j,a_k, a_l\in A$ such that $a_iw$, $a_jw$, $wa_k$, and
$wa_l$ are factors of $W$. Then the $(k+1)$-graph $G_{k+1}$ is
obtained from the follower by means of deleting an arc that
corresponds to one of the four words: $a_iwa_k$, $a_iwa_l$,
$a_jwa_k$, or $a_jwa_l$. Consider the interval which is the
characteristic set for the word, $I_w=[x_w,y_w]$.

Since $w$ is a right special word, we have $I_w\subset
T^{-1}(I_{a_k}\cup I_{a_l})$; since $w$ is a left special word, we have
$I_w\subset T(I_{a_i}\cup I_{a_j})$.

Suppose that point $A\in [0,1]$ partitions $I_w$ into two intervals whose
images lie in $I_{a_k}$ and $I_{a_l}$, respectively, and point
$B\in [0,1]$ partitions it into intervals whose preimages lie in
$I_{a_i}$ and $I_{a_j}$, respectively.

The choice of the minimal non-occurring word (and, hence, the arc to be deleted)
is determined by the mutual location of points $A$ and $B$
and by whether the mapping preserves orientation on these sets or they are reversed.
There are 8 cases, which are divided into four pairs
which correspond to similar sets of words:

\begin{enumerate}

\item $B<A$, $T^{-1}([x_w,B])\subset I_{a_i}$, $T([x_w,A]\subset
I_{a_k} )$

\item $B<A$, $T^{-1}([x_w,B])\subset I_{a_j}$, $T([x_w,A]\subset
I_{a_k} )$

\item $B<A$, $T^{-1}([x_w,B])\subset I_{a_i}$, $T([x_w,A]\subset
I_{a_l} )$

\item $B<A$, $T^{-1}([x_w,B])\subset I_{a_j}$, $T([x_w,A]\subset
I_{a_l} )$

\item $B>A$, $T^{-1}([x_w,B])\subset I_{a_i}$, $T([x_w,A]\subset
I_{a_k} )$

\item $B>A$, $T^{-1}([x_w,B])\subset I_{a_j}$, $T([x_w,A]\subset
I_{a_k} )$

\item $B>A$, $T^{-1}([x_w,B])\subset I_{a_i}$, $T([x_w,A]\subset
I_{a_l} )$

\item $B>A$, $T^{-1}([x_w,B])\subset I_{a_j}$, $T([x_w,A]\subset
I_{a_l} )$

\end{enumerate}

Cases $1$ and $5$ correspond to the forbidden word $a_iwa_k$.

Cases $2$ and $6$ correspond to the forbidden word $a_iwa_l$.

Cases $3$ and $7$ correspond to the forbidden word $a_jwa_k$.

Cases $4$ and $8$ correspond to the forbidden word $a_jwa_l$.

Two pairs of cases correspond to simultaneous reversal or
preservation of orientation of the mapping on characteristic
intervals; two other pairs correspond to opposite orientations.

In the case where the interval exchange transformation preserves
orientation, we have only two possibilities for the arc to be
deleted.

If the mapping is allowed to reverse the intervals
in the process of the interval exchange transformation, all four cases are possible.

Let us introduce the notion of the {\it labeled Rauzy graph}. A
Rauzy graph is said to be {\it labeled} if

\begin{enumerate}
\item The arcs of any crotch are assigned symbols $l$ (``left'') and $r$ (``right'');

\item Some vertices are assigned symbol ``--''.
\end{enumerate}

The {\it follower} of the labeled Rauzy graph is the directed graph
which is the follower of the latter (considered a Rauzy graph with
the labeling neglected) and whose arcs are labeled according the
following rule:
\begin{enumerate}

\item Arcs that enter a crotch should be labeled by the same
symbols as the arcs that enter any left successor of this vertex;

\item Arcs that go out of a crotch should be labeled by the same
symbols as the arcs that go out of any right successor of this
vertex;

\item If a vertex is labeled by symbol ``--'', then all its right
successors should also be labeled by symbol ``--''.

\end{enumerate}

{\bf Remark.} Now let us explain the meaning of the labeled Rauzy
graph. Let the arcs of the incoming crotch correspond to $a_i$ and
$a_j$ and symbols $l$ and $r$ correspond to the left and right set
in the pair $(T(I_{a_i}),T(I_{a_j}))$. If symbols $a_k$ and $a_l$
correspond to the arcs of the outgoing crotch, then symbols $l$
and $r$ appear in accordance with the ``left--right'' order in the
pair $(I_{a_k},I_{a_l})$. A vertex is assigned symbol ``--'' if
the characteristic set that corresponds to it belongs to an
interval which is reversed in the process of interval exchange
transformation.

Below, we give a condition for passing from graph $G_n$ to $G_{n+1}$.

\begin{proposition}
\begin{enumerate}
\item If the graph contains no double crotches that correspond to bispecial factors,
then, when passing from $G_n$ to $G_{n+1}$, we have $G_{n+1}=D(G_n)$.

\item If the vertex that corresponds to a bispecial word is not
labeled by symbol ``--'', then the arcs that correspond to the
forbidden words are chosen among the pairs $lr$ and $rl$.

\item If a vertex is labeled by symbol ``--'', then the arcs to
be deleted should be chosen among the pair $ll$ or $rr$.
\end{enumerate}
\end{proposition}

The evolution of labeled Rauzy graphs is said to be {\it correct}
if {\bf Rules 1} and {\bf 2} are complied with by all graphs in
the evolution starting from $G_1$; the evolution is said to be
{\it asymptotically correct} of {\bf Rules 1} and {\bf 2} are
complied with starting from a certain $G_n$. We say that the
evolution of labeled Rauzy graphs is {\it oriented} if the
$k$-graphs contain no vertices labeled by symbol ``--''.

The definition of the asymptotically correct evolution of Rauzy graphs
allows us to formulate the conditions for a word to be generated by
interval exchange transformation.

\begin{proposition}\label{Osn} A uniformly recurrent word $W$
\begin{enumerate}
\item is generated by interval exchange transformation if the word
is provided with the asymptotically correct evolution of labeled
Rauzy graphs; \item is generated by orientation-preserving
interval exchange transformation if the word is provided with the
asymptotically correct oriented evolution of labeled Rauzy
graphs.
\end{enumerate}
\end{proposition}

Our main result consists in replacing ``if'' with ``if and only if''
in the proposition formulated above.

\subsection{Construction of the dynamical system}
Let us demonstrate that conditions of Theorem~\ref{Osn} are
sufficient for the word to be generated by interval exchange
transformation. First, we show that the word $W$ that satisfies the
hypotheses of Theorem~\ref{Osn} can be the evolution of a certain
point under the following piecewise-continuous transformation of the
segment $T: I\to I$:

\begin{enumerate}
\item $I=[x_0,x_1]\cup [x_1,x_2] \cup \ldots \cup [x_{n-1},x_n]$,
$x_0=0$, $x_n=1$.

\item $I=[y_0,y_1]\cup [y_1,y_2] \cup \ldots \cup [y_{n-1},y_n]$,
$y_0=0$, $y_n=1$.

\item $\sigma \in {S_n}$ is a permutation of a set of $n$ elements.

\item Transformation $T$ maps $(x_i,x_{i+1})$
onto $(y_{\sigma(i)},y_{\sigma(i)+1})$ in a continuous and bijective manner.

\end{enumerate}

Then we demonstrate that the case of piecewise-continuous transformation
can be reduced to interval exchange transformation.

WE construct the piecewise-continuous transformation $T$ step by step.
at the first iteration, we partition the interval into arbitrary
subintervals which correspond to appropriate symbols.
To construct the mapping, it is sufficient to determine the trajectory
of these points and then, for reasons of recurrence, extend it by continuity
to the entire interval.

{\bf Notation.} By virtue of continuity and bijectivity,
the mappings on the intervals are monotonic functions. We shall consider two cases:

\begin{enumerate}

\item All mappings of the intervals are increasing functions.
Such transformation is said to be {\it orientation-preserving}.

\item There are both increasing and decreasing mappings of the intervals.
In this case, we say that transformation {\it does not preserve orientation}.
\end{enumerate}

Let us partition the segment $I=[0,1]$ into $n=\Card A$ arbitrary intervals,
which will be regarded as characteristic sets for the symbols of alphabet $A$:
$[0,1]=I_{a_1}\cup I_{a_2}\cup \ldots \cup I_{a_n}$.

The correspondence between the intervals and the symbols
is defined by Proposition~\ref{Sosed} below.

Let us assume that the mapping $T$ is continuous on each set
$I_{a_i}$, i.e., mapping $T$ can be discontinuous only at the
endpoints of the characteristic sets. The intervals of the
characteristic sets (or their images) that have a common point are
said to be {\it adjacent}.

{\bf Remark.} We can always enlarge the alphabet in such a way that
characteristic sets in the extended alphabet be organized exactly as described above.
The graphs of the $k$-words in the original alphabet then correspond to
$1$-graphs of the extended alphabet and the evolutions of graphs
coincide starting from this moment.

One can directly verify the following assertion.

\begin{proposition}
The images of two adjacent sets under transformations $T$ and $T^{(-1)}$ either are adjacent,
or cannot cover the entire interval.
\end{proposition}

\begin{proposition}\label{Sosed}
The partition intervals can be put in correspondence with the
symbols of the alphabet in such a way that, if $w$ is a special
right $1$-word and $wa_i$, $wa_j$ are subwords, then sets
$I_{a_i}$ and $I_{a_j}$ are adjacent. Similarly, if $w$ is a
special left $1$-word and $a_kw$, $a_lw$ are subwords, then sets
$I_{a_k}$ and $I_{a_l}$ are also adjacent.
\end{proposition}

\Proof Consider an arbitrary crotch in the $1$-graph. Suppose that
arcs that go out of this crotch lead to vertices that correspond to symbols $a_i$ and $a_j$.
Then we set $I_{a_i}=[x_0,x_1]$, $I_{a_j}=[x_1,x_2]$,

Since the characteristic sets are intervals, there is only one
order relation that can be introduced on them (namely, $I_{a_i}< I_{a_j}$ if $x_i<x_j$);
the same order relation can be introduced on their images.
If a pair of characteristic sets reverse their order under transformation $T$,
then we say that, on this pair, the transformation changes orientation.

Consider the images of the intervals under the mapping $T^{-1}$. It
is clear that, if symbol $a_i$ is not a right special $1$-word and
is inevitably followed by symbol $a_j$, then $I_{a_i}\subset
T^{-1}(I_{a_j})$; if it is not a left special $1$-word and is always
preceded by symbol $a_k$, then $T^{-1}(I_{a_i})\subset I_{a_k}$.

In the case where symbol $a_i$ is a right special $1$-word and can be followed
by symbols $a_j$ and $a_k$, we have $I_{a_i}\subset
T^{-1}(I_{a_j}\cup I_{a_k})$; if it is a left special one, then
$T^{-1}(I_{a_i})\subset I_{a_j}\cup I_{a_k}$.

Denote the set of the images of the interval ends under the mapping $T^{-n}$ by $I^n$
and denote the set of the ends of the intervals of characteristic sets by $I^{0}$, i.e., $I^{0}=\{x_0,x_1,\ldots,
x_n\}$.

As follows from considerations in Section~\ref{Sootv:2}, the set that corresponds
to the word $w=w_1w_2\cdots w_n$ is $I_w=I_{w_n}\cap T^{-1}(I_{w_{n-1}})\cap \ldots \cap T^{(-n+1)}(I_{w_1})$;
accordingly, the set of boundary points of the sets that correspond to words of length $n$
is $I^{0}\cup I^{1}\cup \ldots \cup I^{(n-1)}$.

If a right special word is not a bispecial one
(which means that it is not at the same time a left special one),
then the location of the point that partitions this characteristic set
does not matter and it can be chosen arbitrarily.

In the case of orientation-preserving transformation,
the partition should be in agreement with orientation-preserving rules.

In the case where transformation does not preserve orientation,
it is necessary that, in the process of evolution,
the number of ``breaks'' inside the intervals being mapped be finite.

Thus, we can define the images of points on a certain subset
$N\subset I$. It follows from the construction that there exist
intervals $I_k=(x_k,x_{k+1})$ such that,
inside these intervals, our transformation is monotonic.
We can always extend it by continuity to a mapping of the interval into itself.
The resulting piecewise-continuous transformation
is just what we looked for. Let us denote it by $T$.
Note that the initial point whose evolution is the desired word
$W=\{w_n\}$ is the point of intersection for the sequence of
nested intervals which correspond to prefixes $w_0, w_0w_1, w_0w_1w_2, \ldots$.

Thus, we have proved

\begin{theorem}
For a recurrent word $W$ to be generated by
pi\-ec\-e\-wi\-se--con\-ti\-nu\-o\-us
or\-i\-en\-ta\-ti\-on--pre\-ser\-ving transformation, it is
necessary and sufficient that the word be provided with
asymptotically correct oriented evolution of labeled Rauzy
graphs.
\end{theorem}

In the case of orientation-changing transformation, we have

\begin{theorem}
For a recurrent word $W$ to be generated by
pi\-ec\-e\-wi\-se--con\-ti\-nu\-o\-us transformation, it is
necessary and sufficient that the word be provided with
asymptotically correct evolution of labeled Rauzy graphs.
\end{theorem}

\section{Equivalence of the set of uniformly recurrent words generated
by pi\-e\-ce\-wi\-se--con\-ti\-nu\-o\-us transformation to the set
of words generated by interval exchange transformation}

First, let us pass to the dynamics in which almost all points (in
the sense of the Lebesgue measure) have distinct essential
evolutions. Consider the essential evolutions of the points under
transformation $T$.

Consider a piecewise-continuous transformation of intervals. The
theorem of existence of invariant measure (see~\cite{Sin}) yields
that any mapping of a compact set has invariant probabilistic
measure. Therefore, the mapping $T$ has invariant measure $\mu$ and
we can introduce a new semimetric $d(x_1,x_2)=\mu((x_1,x_2))$ on the
segment. Note that the inequality $\mu(x_1,x_2)>0$ does not imply
that the points have different essential evolutions, since the
mapping constructed can be discontinuous. Moreover we need chose a
suitable measure.

The partition $U_1, \ldots, U_n$ of interval is called a {\it pure
partition} if following conditions hold: 1) for any $a_i$
characteristic set $U_i$ is convex, i.e. $U_i$ is closed,
semienclosed or open interval. 2) if two points $x_1,x_2$ has a same
color and the interval $(x_1,x_2)$ contain a break point then images
$T(x_1)$ and $T(x_2)$ has a different colors.

Let $U_1,U_2,\ldots,U_n$ be a partition into characteristic sets.
The partition $V_1, V_2,\ldots, V_m$ is called {\it subpartition} if
each characteristic set $U_i$ is union of sets $V_j$:
$U_i=V_{i_1}\cup V_{i_2}\cup V_{i_j}$.

Let $W$ be an evolution for some partition and $W'=\{w'_i\}$ be an
evolution for its subpartition. It is clear that $W$ raise from $W'$
by gluing of symbols.  It is easy to see that every partition has a
pure subpartition.

Let $W$ would be an uniformly recurrent word, corresponding some
evolution with partition $\cal U$, $\widehat{W}$ is a word,
corresponding the evolution with the same initial point and
subpartition ${\cal U}'$. Gluing morphism of alphabets give us a
natural morphism $\pi$ of the words $\pi: \widehat{W}\to W$.
Logically, $\widehat{W}$ may not necessary uniformly recurrent, but
there is an u.r. $\widehat{U}$ such that
$\widehat{U}\preceq\widehat{W}$. Then $U=\pi(\widehat{U})\preceq
W=\pi(\widehat{W})$ and hence u.r via theorem \ref{RecCont}. In the
corresponding to $\widehat{U}$ closed orbit would be point
corresponding an evolution with projection $W$. Indeed, let $w$ be
an arbitrary subword of $W$, it occurs in any subword of $W$ of
length $\ge k(w)$. Hence in any subword of $\widehat{U}$ of length
$\ge k(w)$ there exist a subword $\hat{w}$ such that
$\pi(\hat{w})=w$. Hence for any surrounding of zero position of $W$
there exist $\widehat{V}\equiv \widehat{V}$ such that
$\pi(\widehat{V})$ has the same surrounding.  It remains to use
compactness argument. Corresponding word $\widehat{W}'$ will be u.r.
In the sequel we consider only pure partitions.

\begin{proposition}
Suppose that points $x_1<x_2$ have the same evolution and
and the partition is pure. Then any point of the interval
$(x_1,x_2)$ has the same essential evolution.
\end{proposition}

\Proof The fact that$x_1$ and $x_2$ have the same evolution yields
that $T(x_1)$ and $T(x_2)$ have the same evolution as well and the
image of the interval $(x_1,x_2)$ is the interval
$(T(x_1),T(x_2))$.\Endproof

Let $W$ be an nonperiodic uniformly recurrent word generated by a
pi\-e\-ce\-wi\-se--con\-ti\-nu\-o\-us transformation of the
interval. Consider the space of words with shift operator. For $W$
there exist minimal invariant set $N_W$ (see \cite{BBL}). Every
point $V\in N_W$ correspond set of points of system of intervals
with corresponding essential evolution. On the space $N_W$ there
exist invariant (for shift operator) probabilistic measure $\nu$
(\cite{Sin}). This measure induce on interval system invariant
probabilistic measure $\mu$ by the natural way. The more detailed
description of this measure $\mu$ is follows. If $I$ is an interval,
corresponding to the word $u$, and $N$ is the set of points of
closed orbit of $W$ having subword $u$ on given position, then
$\mu(I)=\nu(N)$. If $I$ has no intersection with any such intervals,
then $\mu(I)=0$.

The measure $\mu$ is invariant and induces a semimetric. (Length of
any segment is equal its measure.)

Consider topology that glue point with zero distance and construct
corresponding factor-dynamics. Obtaining map is piecewise-isometric,
i.e it is interval exchange transformations. Every such glued
interval is contained in maximal glued interval that the number of
such maximal intervals is countable and their common measure is
zero. Hence for almost every point (in the sense of $\mu$) of our
compact $M$ has orbit without intersection of any such glued
interval. We have constructed piecewise continuous transformation of
system $M'$ generating an u.r. word, $W'$ which is equivalent to
$W$, and hence, via compactness argument, $M'$ has a point whose
essential evolution is $W$, because if there is a word with
essential evolution $W'$, then any equivalent word can be observed.

We have proved

\begin{theorem}
Suppose that the word $W$ is generated by a
pi\-e\-ce\-wi\-se--con\-ti\-nu\-o\-us transformation of the
interval. Then there exists a word $W'$ which is equivalent to the
given one and is generated by interval exchange transformation.
\end{theorem}

By theorem, this interval exchange transformation involves {\it all}
words that are equivalent to $W'$ including the word $W$ (in the
sense of the essential evolution). Thus, we have proved the
following assertion.

\begin{theorem}
For a recurrent word $W$ to be generated by a
pi\-e\-ce\-wi\-se--con\-ti\-nu\-o\-us
or\-i\-en\-ta\-ti\-on--pre\-ser\-ving transformation, it is
necessary and sufficient that the word be provided with
asymptotically correct oriented evolution of labeled Rauzy graphs.
\end{theorem}

In the case of orientation-changing transformation, we have proved
the following.

\begin{theorem}
For a recurrent word $W$ to be generated by a
pi\-e\-ce\-wi\-se--con\-ti\-nu\-o\-us
or\-i\-en\-ta\-ti\-on--chang\-ing transformation, it is necessary
and sufficient that the word be provided with asymptotically correct
evolution of labeled Rauzy graphs.
\end{theorem}

\end{document}